\documentclass[12pt]{article}

\usepackage{amsthm}
\usepackage{amsmath}
\usepackage{mathrsfs}
\usepackage{stmaryrd}
\usepackage{amssymb}
\usepackage{diagbox}
\usepackage{tikz-cd}
\usepackage{amsfonts}
\usepackage{mathtools}
\usepackage{authblk}
\usepackage{bm}
\usepackage[utf8]{inputenc}
\usepackage{indentfirst}
\usepackage{graphicx}
\usepackage{enumerate}
\usepackage[top=2.54cm, bottom=2.54cm, left=3.28cm, right=3.28cm]{
geometry}
\usepackage{verbatim}
\usepackage[colorlinks,linkcolor=blue]{hyperref}

\newcommand{\be}{e}

\newcommand{\inv}{^{-1}}

\newcommand{\mbf}{\mathbb{F}}

\newcommand{\vp}{\varphi}

\newtheorem{thm}{Theorem}
\theoremstyle{definition}
\newtheorem{defn}[thm]{Definition}
\newtheorem{lem}{Lemma}

\newtheorem{cor}{Corollary}
\newtheorem{exm}{Example}

\newtheorem{remark}{Remark}

\makeatletter
\newcommand{\proofpart}[2]{
  \par
  \addvspace{\medskipamount}
  \noindent\textbf{Part #1. #2}
  \par\nobreak\smallskip
  \@afterheading
}
\makeatother

\begin{document}

\title{Annihilating polynomial, Jordan canonical from, and generalized spectral characterizations of Eulerian graphs }
\author{\small Kunyue Li$^{{\rm a}}$\quad\quad Wei Wang$^{\rm a}$\thanks{Corresponding author: wang\_weiw@xjtu.edu.cn}\quad\quad Hao Zhang$^{\rm b}$
\\
{\footnotesize$^{\rm a}$School of Mathematics and Statistics, Xi'an Jiaotong University, Xi'an 710049, P. R. China}\\
{\footnotesize$^{\rm b}$School of Mathematics, Hunan University, Changsha 410082, P. R. China}}

\date{}
	\maketitle
	\begin{abstract}
Let $G$ be an Eulerian graph on $n$ vertices with adjacency matrix $A$ and characteristic polynomial $\phi(x)$. We show that when $n$ is even (resp. odd), the square-root of $\phi(x)$ (resp. $x\phi(x)$) is an annihilating polynomial of $A$, over $\mathbb{F}_2$. The result was achieved by applying the Jordan canonical form of $A$ over the algebraic closure $\bar{\mathbb{F}}_2$. Based on this, we show a family of Eulerian graphs are determined by their generalized spectrum among all Eulerian graphs, which significantly simplifies and strengthens the previous result.

	\end{abstract}

\section{Introduction}
The spectrum of a graph $G$ is the multiset of the eigenvalues of its adjacency matrix $A(G)$. It is well-known that the spectrum of a graph encodes a lot of combinatorial information
about the given graph, and the interplay between linear algebraic and the combinatorics constitutes the main theme of spectral graph theory; see e.g.,~\cite{BH,CDS}.

A long-standing open question in this area is: ‘‘Which graphs are determined by their spectra (DS for short)?" Here a graph $G$ is DS means that no other graph non-isomorphic to $G$ can share the same spectrum as $G$, e.g., the path $P_n$, the complete graph $K_{n}$, and the regular complete bipartite graph $K_{m,m}$ are some simple examples of DS graphs.
However, it turns out that showing a given graph to be DS is generally very hard and challenging. Despite many efforts, only very few graphs are known to be DS, and the methods involved to show them to be DS heavily depend on the particular structures of these graphs.  
 We refer the reader to \cite{DH1,DH2} for more background and known results on this topic.

 Wang and Xu~\cite{WX1} initiated the study the generalized spectral characterization of graphs. A graph $G$ is said to be \emph{determined by its generalized spectrum} (DGS for short) if whenever $H$ is a graph such that $H$ and $G$ are cospectral with cospectral complements, then $H$ is isomorphic to $G$.

 Let $G$ be a graph with adjacency matrix $A(G)$ on $n$ vertices. Let $W=W(G):=[e, Ae, \ldots, A^{n-1}e]$ be the \emph{walk-matrix} of graph $G$, where $e=(1, 1, \dots, 1)^{\rm T}$ is the all-one vector. $G$ is \emph{controllable} if its walk-matrix $W(G)$ is non-singular.

 In Wang~\cite{Wang2}, the author proved the following theorem.

\begin{thm}[Wang~\cite{Wang2}]\label{wm}
If ${2^{-\lfloor n/2 \rfloor}}{\det W }$ (which is always an integer) is odd and square-free, then $G$ is DGS.
 \end{thm}

 However, the above theorem fails for some interesting family of graphs such as Eulerian graphs. This is because for an Eulerian graph, all the entries (except for the ones in the first column) of $W$ are even, and hence $2^{n-1}$ always divides $\det W$ and
 ${2^{-\lfloor n/2 \rfloor}}{\det W }$ is neither odd nor square-free.

  It is tempting to give an analogous condition to Theorem~\ref{wm} for an Eulerian graph to be DGS.
 Nevertheless, this is generally not available; see Example 1 below. So we are merely satisfied with a slightly weaker result described below.

 Denoted by $\mathcal{E}_n$ the set of all Eulerian graphs on $n$ vertices. An Eulerian graph $G\in{\mathcal{E}_n}$ is DGS \emph{among all Eulerian graphs}, if for any $H\in{\mathcal{E}_n}$, $H$ and $G$ are cospectral with cospectral complements imply that $H$ is isomorphic to $G$.

 The main result of the paper is the following theorem.

 \begin{thm}\label{main}Let $G$ be an Eulerian graph. Then $2^{-\lfloor (3n-3)/2\rfloor}$ always divides ${\det W}$. Moreover, If ${2^{-\lfloor (3n-3)/2\rfloor}}{\det W}$ is odd and square-free, then $G$ is DGS among all Eulerian graphs.
\end{thm}

Let us give an example for illustration.

\begin{exm}
 Let $G$ and $H$ be two graphs with adjacency matrices $A$ and $B$ given as follows, respectively:
$$ A=\left(\begin{array}{cccccccccc}
 0 & 0 & 0 & 0 & 1 & 1 & 1 & 1 & 1 & 1 \\
 0 & 0 & 0 & 0 & 0 & 1 & 0 & 1 & 0 & 0 \\
 0 & 0 & 0 & 1 & 0 & 0 & 0 & 1 & 0 & 0 \\
 0 & 0 & 1 & 0 & 0 & 1 & 0 & 0 & 1 & 1 \\
 1 & 0 & 0 & 0 & 0 & 0 & 0 & 0 & 1 & 0 \\
 1 & 1 & 0 & 1 & 0 & 0 & 0 & 0 & 0 & 1 \\
 1 & 0 & 0 & 0 & 0 & 0 & 0 & 1 & 0 & 0 \\
 1 & 1 & 1 & 0 & 0 & 0 & 1 & 0 & 0 & 0 \\
 1 & 0 & 0 & 1 & 1 & 0 & 0 & 0 & 0 & 1 \\
 1 & 0 & 0 & 1 & 0 & 1 & 0 & 0 & 1 & 0 \\
\end{array}
\right), B=\left(
\begin{array}{cccccccccc}
 0 & 0 & 1 & 1 & 1 & 1 & 0 & 0 & 1 & 1 \\
 0 & 0 & 1 & 0 & 0 & 1 & 0 & 1 & 0 & 1 \\
 1 & 1 & 0 & 0 & 0 & 0 & 0 & 1 & 1 & 0 \\
 1 & 0 & 0 & 0 & 1 & 0 & 0 & 0 & 0 & 0 \\
 1 & 0 & 0 & 1 & 0 & 0 & 1 & 1 & 0 & 0 \\
 1 & 1 & 0 & 0 & 0 & 0 & 0 & 0 & 1 & 1 \\
 0 & 0 & 0 & 0 & 1 & 0 & 0 & 0 & 0 & 0 \\
 0 & 1 & 1 & 0 & 1 & 0 & 0 & 0 & 0 & 0 \\
 1 & 0 & 1 & 0 & 0 & 1 & 0 & 0 & 0 & 0 \\
 1 & 1 & 0 & 0 & 0 & 1 & 0 & 0 & 0 & 0 \\
\end{array}
\right).$$

It is easy to verify that $G$ is an Eulerian graph and $2^{-13}\det W(G)=-43$, which is odd and square-free. Thus, $G$ is DGS among all Eulerian graphs.

Moreover, it can be easily computed that $\phi(G;x)=\phi(H;x)=-4 + 32 x + 48 x^2 - 88 x^3 - 108 x^4 + 62 x^5 + 73 x^6 - 12 x^7 - 17 x^8 + x^{10}$ and
 $\phi(\bar{G};x)=\phi(\bar{H};x)=3 + 20 x - 142 x^3 - 31 x^4 + 176 x^5 + 67 x^6 - 54 x^7 - 28 x^8 + x^{10}$; i.e., $G$ and $H$ are generalized cospectral. Note that $H$ is not an Eulerian graph. Thus, even if $G$ satisfies the conditions of Theorem~\ref{main}, there does exist a non-Eulerian graph that is generalized cospectral with $G$ and non-isomorphic to $G$. Therefore, in certain sense, Theorem~\ref{main} is the best that we can hope.

\end{exm}
\begin{remark}
 It will be soon clear from the proof of Theorem~\ref{main} that the condition that $G$ is connected is not essential; we can replace the condition that $G$ is Eulerian by just requiring that the degree of every vertex of $G$ is even. Moreover, the conclusion still holds for graphs $G$ with the degree of every vertex being odd; this can be seen easily by considering the complement $\bar{G}$ of $G$ (the degree of $\bar{G}$ are all even in this case).  \\
\end{remark}
In~\cite{QJW1}, Qiu, Ji and the second author considered the generalized spectral characterizations of Eulerian graphs based on the methods developed in \cite{Wang1,Wang2}. Theorem~\ref{main} significantly simplifies and strengthens the result of \cite{QJW1}, the proof of which follows the main strategy introduced in Qiu et al.~\cite{QWWZ} for showing a graph to be DGS.

The new ingredient is a discovery of an interesting theorem (see Theorem~\ref{thm:mainresult}) which asserts that for an Eulerian graph $G$ with adjacency matrix $A$ and characteristic polynomial $\phi(x)$, when $n$ is even (resp. odd), the square-root of $\phi(x)$ (resp. $x\phi(x)$) is an annihilating polynomial of $A$, over $\mathbb{F}_2$. The result is of independent interest and somewhat unexpected, the proof of which was achieved by using Jordan canonical form of $A$ over $\bar{\mathbb{F}}_2$ (the algebraic closure of $\mathbb{F}_2$).

 The rest of the paper is organized as follows. In Section 2, we provide some preliminary results. Section 3 is dedicated to the proof of Theorem~\ref{main}. The proof of the key intermediate result - Theorem~\ref{thm:mainresult}, is presented in Section 4.

 \section{Preliminaries}

 For convenience of the reader, we give some preliminary results that are needed later in the proof Theorem~\ref{main}.

 \subsection{Notations and terminologies}

 Let $G=(V, E)$ be a graph with vertex set $V=\{v_1, v_2, \ldots, v_n\}$ and edge set $E$.
 The \emph{adjacency matrix} of $G$ is an $n\times n$ matrix $A(G)=(a_{ij})$, where $a_{ij}=1$ if $v_i$ and $v_j$ are adjacent, and $a_{ij}=0$ otherwise.
  The \emph{characteristic polynomial of $G$} is defined as $\phi(x)=\phi(G; x)=\det(x I- A(G))$.
  The \emph{adjacency spectrum} of $G$, denoted by ${\rm Spec}(G)$, is the multiset of all the eigenvalues of $A(G)$.
  Two graphs $G$ and $H$ are \emph{cospectral} if ${\rm Spec}(G)={\rm Spec}(H)$.
 A graph $G$ is said to be \emph{determined by the spectrum} (DS for short), if any graph that is cospectral with $G$ is isomorphic to $G$, i.e., ${\rm Spec}(G)={\rm Spec}(H)$ implies that $H$ is isomorphic to $G$ for any $H$.

  Two graphs $G$ and $H$ are \emph{generalized cospectral} if they are cospectral with cospectral complement, i.e., ${\rm Spec}(G)={\rm Spec}(H)$ and ${\rm Spec}(\overline{G})={\rm Spec}(\overline{H}).$
   A graph $G$ is said to be \emph{determined by the generalized spectrum} (DGS for short), if any graph that is generalized cospectral with $G$ is isomorphic to $G$, i.e., $G$ is DGS if ${\rm Spec}(G)={\rm Spec}(H)$ and ${\rm Spec}(\overline{G})={\rm Spec}(\overline{H})$ imply that $H$ is isomorphic to $G$ for any $H$.

  In this paper, we are mainly concerned with \emph{Eulerian graphs} $G$, i.e., $G$ is a graph having an Eulerian trail, or equivalently, $G$ is connected and the degree of every vertex is even.\\

  \noindent
  \textbf{Notations}: Throughout, let $p$ be a prime number. We use ${\rm rank}_p\,M$ to denote the rank of an integral matrix $M$ over the finite field $\mathbb{F}_p$. Also, for two integers $a$ and $b$, we write $a\equiv b~({\rm mod}~ p)$ and $a=b$ over $\mathbb{F}_p$ interchangeably.

 \subsection{Main strategy}

 In this subsection, we shall briefly describe the main strategy for showing a graph to be DGS; see~\cite{QJW1}-\cite{Wang2} for some details.

 First, we give some definitions. An $n$ by $n$ matrix $Q$ is an \emph{orthogonal matrix} if $Q^{\rm T}Q=I_n$, and it is a \emph{rational} orthogonal matrix if all the entries of $Q$ are rational numbers.
 An orthogonal matrix $Q$ is called \emph{regular}, if every row (column) sum of $Q$ is one, i.e., $Qe=e$, where $e$ is the all-one vector.

 The following theorem gives a simple characterization for two graphs to be generalized cospectral.

\begin{thm} [Johnson and Newman~\cite{JN}; Wang and Xu~\cite{WX1}]\label{Ch}
A pair of graphs $G$ and $H$ are generalized cospectral if and only if there exists a regular orthogonal matrix $Q$ such that
\begin{equation*}
Q^{\rm T}A(G)Q=A(H).
  \end{equation*}
  Moreover, if $G$ be a controllable, then $Q=W(G)W(H)^{-1}$ and hence is rational and unique.
\end{thm}

 \begin{defn}
 $\mathcal{Q}(G):=\{Q\in{RO_n(\mathbb{Q})}:\,Q^{\rm T}A(G)Q ~{\rm is~a~(0,1)-matrix}\}$, where $RO_n(\mathbb{Q})$ denotes the set of all regular rational orthogonal matrices of order $n$.
  \end{defn}

  \begin{thm}[Wang and Xu~\cite{WX1}]
  Let $G$ be a controllable graph. Then $G$ is DGS if and only if $\mathcal{Q}(G)$ contains only permutation matrices.
  \end{thm}

  \begin{defn}Wang and Xu~\cite{WX1}Let $Q$ be a regular rational orthogonal matrix. The level of $Q$, denoted by $\ell$ or $\ell(Q)$,  is the smallest positive integer $k$ such that $k Q$ is an integral matrix.
\end{defn}

 For example, let
$$Q_1=
\frac{1}{2}\left[ \begin{array}{cccc}
-1&1&1&1\\
1&-1&1&1\\
1&1&-1&1\\
1&1&1&-1
\end{array}\right], Q_2=\frac{1}{3}\left[ \begin{array}{cccccc}
2&-1&-1&1&1&1\\
-1&2&-1&1&1&1\\
-1&-1&2&1&1&1\\
1&1&1&2&-1&-1\\
1&1&1&-1&2&-1\\
1&1&1&-1&-1&2
\end{array}\right].$$
Then $Q_1$ and $Q_2$ are two regular rational orthogonal matrix of level 2 and 3, respectively.

Clearly, a regular rational orthogonal matrix is a permutation matrix if and only if $\ell=1$. Thus, the main strategy for us to show a given graph $G$ to be DGS is to show
that for every $Q\in {\mathcal{Q}(G)}$, the level of $Q$ is one.
\subsection{Smith normal form and the level of $Q$}
It turns out that the level of $Q\in {\mathcal{Q}(G)}$ is closely related to the $n$-th invariant factor of the walk matrix of $G$, which will be briefly described below.

Recall that an $n \times  n$ matrix $U$ with integer entries is \emph{unimodular} if $\det U = \pm 1$. For every integral matrix $M$ with full rank, there exist two unimodular matrices
$U$ and $V$ such that $M=USV$, where $S={\rm diag}(d_{1}, d_{2}, \dots, d_{n})$ with $d_{i} \mid d_{i+1}$ for $i = 1, 2, \dots , n-1$.  The diagonal matrix $S$ is known as the \emph{Smith Normal Form} (SNF for short) of $M$, and $d_i$ is the $i$-th \emph{invariant factor} of $M$.

\begin{lem}[\cite{QWWZ}]\label{level}
 Let $Q$ be a rational orthogonal matrix with level $\ell$. Suppose $QX=Y$ for two non-singular integral matrices $X$ and $Y$. Then $\ell\mid d_n(X)$,
 where $d_n(X)$ is $n$-th invariant factor of $X$.
 \end{lem}

As a simple consequence, we have the following

 \begin{cor} [\cite{WX1}]\label{L11} Let $G$ be any graph and $Q\in {\mathcal{Q}(G)}$ with level $\ell$. Then $\ell\mid d_n(W)$, where $d_n(W)$ is the last invariant factor of $W$; in particular, $\ell\mid \det W$.
\end{cor}

Let $$\Sigma_n:=\{G\in {\mathcal{E}_n}:\,2^{-\lfloor (3n-3)/2\rfloor}\det W ~{\rm is~ odd~ and~ squarefree}\}.$$

 \begin{lem}[\cite{QJW1}]\label{LL2}
 Let $G\in{\Sigma_n}$. Then the SNF of the walk matrix $W$ of $G$ is as follows:
 $${\rm diag}(\underbrace{1,2,2,\ldots,2}_{\lceil \frac{n+1}2\rceil},\underbrace{4,4,\ldots,4,4b}_{\lfloor \frac{n-1}2\rfloor}),$$
 where $b$ is an odd and square-free integer.
\end{lem}

Let $G\in{\Sigma_n}$ and $Q\in {\mathcal{Q}(G)}$ with level $\ell$. Then according to Lemmas~\ref{L11} and \ref{LL2}, we have $\ell \mid 4b$. Moreover, we have the following

\begin{lem}[\cite{Wang1}]\label{LLL1}  Let $p$ be an odd prime. For any graph $G$, suppose that $p\mid \det W$ and $p^2\nmid \det W$. Then $p\nmid \ell$.
\end{lem}

Thus by Lemma~\ref{LLL1}, any odd prime $p$ is not a divisor of $\ell$. Thus, we further have $\ell \mid 4$. We record this as the following

\begin{lem}\label{Le1}  Let $G\in{\Sigma_n}$ and $Q\in {\mathcal{Q}(G)}$ with level $\ell$. Then $\ell \mid 4$.
\end{lem}

So in the remaining, we shall show that $\ell$ has to be odd, and as a result, $\ell=1$ and $G$ is DGS.

\section{Proof of Theorem~\ref{main}}

In this section, we present the proof of Theorem~\ref{main}. First, we need some preparations.



 \begin{defn}
Define $$\bar{W}:=[e,\frac{Ae}2,\ldots,\frac{A^{n-1}e}2].$$
\end{defn}
Note that $Ae\equiv0~({\rm mod}~2)$. Thus, $\bar{W}$ is an integral matrix.

\begin{lem}[Lemma 18 in \cite{QJW1}] \label{PQ} Let $G\in{\Sigma_n}$. Then the SNF of $\bar W$ is as follows:
$${\rm diag}(\underbrace{1,1,\ldots,1}_{\lceil\frac{n+1}2\rceil},\underbrace{2,2,\ldots,2,2b}_{\lfloor\frac{n-1}2\rfloor}).$$

\end{lem}

\begin{lem}[c.f.~\cite{QJW1}] Let $r=\lceil\frac{n+1}2\rceil$. Then the first $r$ columns of $\bar{W}$ are linearly independent over $\mathbb{F}_2$.

\end{lem}
\begin{proof} By Lemma~\ref{PQ}, we have ${\rm rank}_2\,\bar{W}=r$. For contradiction, suppose that the first $r$ columns of $\bar{W}$ are linearly dependent over $\mathbb{F}_2$.
Let $s~(1\leq s<r)$ be the maximum integer such that the first $s$ columns of $\bar{W}$ are linearly independent over $\mathbb{F}_2$.
Write $\beta=\frac{Ae}2$ and $\bar{W}=[e,\beta,A\beta,\ldots,A^{n-2}\beta]$. Then $A^s\beta\in {\rm span}<e, \beta,A\beta,\ldots,A^{s-1}\beta>$.
Thus $$A^s\beta=c_{-1}e+c_0\beta+c_1A\beta+\cdots+c_{s-1}A^{s-1}\beta,$$
for some $c_{-1},c_0,\ldots,c_{s-1}\in {\mathbb{F}_2}$. It follows that $$A^{s+1}\beta=c_{-1}Ae+c_0A\beta+c_1A^2\beta+\cdots+c_{s-1}A^{s}\beta\in {\rm span}<e, \beta,A\beta,\ldots,A^{s-1}\beta>.$$ Similarly, we can show $$A^{s+t}\beta\in {\rm span}<e, \beta,A\beta,\ldots,A^{s-1}\beta>$$ for any $t\geq 2$.
Thus, ${\rm rank}_p \,\bar{W}=s<r$; a contradiction. This completes the proof.

\end{proof}



 Recall that an \emph{elementary subgraph} of a graph $G$ with $i$ vertices is a disjoint union of cycles and edges on a total of $i$ vertices. Denote the set of elementary subgraphs of $G$ with $i$ vertices by $\mathcal{H}_{i}$.

 \begin{lem}[Sachs' Coefficients Theorem]
Let $\phi(x)=x^{n}+c_{1}x^{n-1}+\cdots+c_{n-1}x+c_{n}$ be the characteristic polynomial of the graph $G$, and let $\mathcal{H}_{i}$ be the set of elementary subgraphs of $G$ with $i$ vertices. Then
$$c_{i}=\sum_{H\in\mathcal{H}_{i}}(-1)^{p(H)}2^{c(H)}~(i=1,\ldots,n),$$
where $p(H)$ denotes the number of components of $H$ and $c(H)$ denotes the number of cycles in $H$.
\end{lem}

First suppose that $n$ is even. By Sachs' Coefficient Theorem, $c_{i}$ is even when $i$ is odd. Thus,
 \begin{eqnarray*}
   \phi(x) &\equiv & x^{n}+c_{2}x^{n-2}+\cdots+c_{n-2}x^{2}+c_{n} \\
        &\equiv& (x^{n/2}+c_2x^{n/2-1}+\cdots+c_{n-2}x+c_n)^2\pmod2.
     \end{eqnarray*}
 Similarly, if $n$ is odd, then $$\phi(x) \equiv  x(x^{(n-1)/2}+c_{2}x^{(n-3)/2}+\cdots+c_{n-3}x+c_{n-1})^2 \pmod2.$$
Define
\begin{equation}\label{PQW}
\vp(x):=\begin{cases}
&x^{\frac{n}{2}}+c_{2}x^{\frac{n-2}{2}}+\cdots+c_{n-2}x+c_{n}, ~\text{if}~ n~\text{is~even};\\
&x^{\frac{n+1}{2}}+c_{2}x^{\frac{n-1}{2}}+\cdots+c_{n-3}x^2+c_{n-1}x, ~\text{if}~n~\text{is~odd}.
\end{cases}
\end{equation}
\begin{thm}[\cite{Wang2}] Let $\vp$ be defined as the above. Then we have $\phi=\vp^2$ or $x\phi=\vp^2$ according to whether $n$ is even or odd.
\end{thm}

The following theorem shows that $\vp$ is always the annihilating polynomial of the adjacency matrix $A$ of an Eulerian graph $G$ over $\mathbb{F}_2$, which lies at the heart of the proof of Theorem~\ref{main}.

\begin{thm}\label{thm:mainresult}
    Let $A$ be an adjacency matrix of an Eulerian graph $G$. Let $\phi$ be the characteristic polynomial of $A$ and $\vp$ be defined as in Eq.~\eqref{PQW}. Then we have $\vp(A)=0$ over $\mbf_2$.
\end{thm}

We shall postpone the proof of Theorem~\ref{thm:mainresult} to the next section. Assuming at the moment that Theorem~\ref{thm:mainresult} holds, we present the proof of Theorem~\ref{main}.

 \begin{proof}[Proof of Theorem~\ref{main}]
 We only prove the case that $n$ is even, the case that $n$ is odd can be proved in a similar way (see Remark 2 below).

 Let $H$ be any Eulerian graph that is generalized cospectral with $G$. We shall show $H$ is isomorphic to $G$. By Theorem~\ref{Ch}, there exists a regular rational orthogonal matrix
 $Q\in {\mathcal{Q}(G)}$ with level $\ell$ such that $Q^{\rm T}AQ=B$, where $A$ and $B$ are the adjacency matrices of $G$ and $H$, respectively.
By Lemma~\ref{Le1}, we known that $\ell$ divides 4. It remains to show that $\ell$ is odd.

 It follows from $Q^{\rm T}AQ=B$ and $Qe=e$ that $Q^{\rm T}A^ke=B^ke$ for any $k\geq 0$. Thus we have
  \begin{gather} \nonumber
Q^{\rm T}e=e,\\\nonumber
Q^{\rm T}\frac{Ae}2=\frac{Be}2,\\\nonumber
\vdots\\
Q^{\rm T}\frac{A^{n/2}e}2=\frac{B^{n/2}e}2,\\ \nonumber
Q^{\rm T}\frac{A^{n/2+1}e}2=\frac{B^{n/2+1}e}2,\\\nonumber
\vdots\\\nonumber
Q^{\rm T}\frac{A^{n-1}e}2=\frac{B^{n-1}e}2.\nonumber
\end{gather}

 Let $\vp(x)=x^{n/2}+c_2x^{n/2-1}+\cdots+c_{n-2}x+c_n$. By Theorem~\ref{thm:mainresult} we get that $\vp(A)=0$ over $\mathbb{F}_2$. Thus we have
 \begin{equation}
 \vp(A)Ae\equiv 0~({\rm mod}~4),
 \end{equation}
 or equivalently,
\begin{eqnarray}
 c_n\frac{Ae}2+c_{n-2}\frac{A^2e}2+\cdots+c_2\frac{A^{n/2}e}2+\frac{A^{n/2+1}e}2\equiv 0~({\rm mod}~2).
\end{eqnarray}

Since $G$ and $H$ are cospectral, we have $\vp(B)=0$ over $\mathbb{F}_2$. Similarly, we have
\begin{eqnarray}
 c_n\frac{Be}2+c_{n-2}\frac{B^2e}2+\cdots+c_2\frac{B^{n/2}e}2+\frac{B^{n/2+1}e}2\equiv 0~({\rm mod}~2).
\end{eqnarray}

 Let $r:={\rm rank}_2\,W=\lceil \frac{n+1}2\rceil=n/2+1$. Multiplying both sides of the second, the third,$\ldots$, the $r$-th equalities in Eq.~(2) by $c_{n}/2,c_{n-2}/2,\ldots,c_2/2$, respectively, then adding them to the $(r+1)$-th equality gives
 that $Q^{\rm T}\frac{\vp(A)Ae}4=\frac{\vp(B)Be}4$. Similarly, we have $Q^{\rm T}\frac{\vp(A)A^ie}4=\frac{\vp(B)B^ie}4$ for $i=1,2,\ldots,n/2-1$, i.e.,

 $$Q^{\rm T}\hat{W}(G)=\hat{W}(H),$$
 where $$\hat{W}(G):=[e,\frac{Ae}2,\ldots,\frac{A^{n/2}e}2,\frac{\vp(A)Ae}4,\ldots,\frac{\vp(A)A^{n/2-1}e}4]$$ and
 $$\hat{W}(H):=[e,\frac{Be}2,\ldots,\frac{B^{n/2}e}2,\frac{\vp(B)Be}4,\ldots,\frac{\vp(B)B^{n/2-1}e}4].$$

 Note that $\det \hat{W}(G)=2^{-(\frac{3n}2-2)}\det W(G)$, which is odd. Further note that both $\hat{W}(G)$ and $\hat{W}(H)$ are integral matrices.
 It follows from Lemma~\ref{level} that $\ell\mid \det \hat{W}(G)$. Therefore $\ell$ is odd and $\ell=1$ and $Q$ is a permutation matrix. This shows $H$ is isomorphic to $G$.
 The proof is complete.
 \end{proof}

 \begin{remark} When $n$ is odd, the proof of Theorem~\ref{main} follows the main lines of the case that $n$ is even, except that we need some slight modifications to Eq.~(3).

 Let $\phi=x^{k}\phi_1^2(x)$ over $\mathbb{F}_2$, where $\gcd(x,\phi_1(x))=1$ and $k\geq 1$ is odd. If $k=1$, then it follows from Theorem~\ref{thm:mainresult} that $\varphi(A)=A\phi_1(A)=0$. Note that ${\rm rank}_2\,A=n-1$ and $\phi_1(A)$ is symmetric. Then we have $\phi_1(A)=J_n$ (the all-one matrix of order $n$).
 Thus, $$\varphi(A)e=A\phi_1(A)e=\phi_1(A)(Ae)=(J_n+2M)Ae\equiv(e^{\rm T}Ae)e~({\rm mod}~4),$$ where $M$ is some integral matrix. If $k\geq 3$, then we have $A^{\frac{k-1}2}\phi_1(A)=0$ over $\mathbb{F}_2$ (see Theorem~\ref{odd} in Section 4). It follows that $$\varphi(A)e=A (A^{\frac{k-1}2}\phi_1(A))e=(A^{\frac{k-1}2}\phi_1(A))(Ae)\equiv0~({\rm mod}~4).$$

 In summary, if $k=1$, let $\tilde{\varphi}(x)=\varphi(x)-e^{\rm T}Ae$, and $\tilde{\varphi}(x)=\varphi(x)$ if $k\geq 3$. Then $\hat{\varphi}(A)e\equiv 0~({\rm mod}~4)$.
 Since $H$ is cospectral with $G$, we have $e^{\rm T}Ae=e^{\rm T}Be$ and hence $\hat{\varphi}(B)e\equiv 0~({\rm mod}~4)$ still holds.
 \end{remark}

\section{Proof of Theorem~\ref{thm:mainresult}}

In this section, we present the proof of Theorem~\ref{thm:mainresult}. First, we shall focus on the case that $n$ is even (assume henceforth that $n$ is even, unless stated otherwise). 

Before giving the proof of Theorem~\ref{thm:mainresult}, we fix some notations. Let $\mbf$ be the splitting field of $\phi$ over $\mbf_2$. Suppose $J$ is the Jordan normal form of $A$ over $\mbf$ with $P\in M_{n\times n}(\mbf)$ such that $A=PJP\inv$. Suppose
\[J=\begin{pmatrix}
    J_1&&\\
    &\ddots&\\
    &&J_t
\end{pmatrix},\]
where $J_i=J(\lambda_i,\ell_i)$ is the Jordan block of order $\ell_i$. Suppose
\[P=(\alpha_{11},\dots,\alpha_{1\ell_1},\dots,\alpha_{t1},\dots,\alpha_{t\ell_t}).\]
Since $A\be=0$, without loss of generality, we may assume that $\alpha_{11}=\be$.

Let $Q=P^{\rm T}P$, then from $A^{\rm T}=A$, we see that $J^{\rm T}Q=QJ$. This implies that $Q$ and $QJ$ are both symmetric and $Q$ is invertible. Suppose
\[Q=\begin{pmatrix}
    Q_{11}&\cdots&Q_{1t}\\
    \vdots&&\vdots\\
    Q_{t1}&\cdots&Q_{tt}
\end{pmatrix},\]
has the same partition as $J$. Now we prove some useful lemmas.

\begin{lem}\label{lem:ealpha0}
    For any $\alpha\in \mbf^n$, $\alpha^{\rm T}\alpha=0$ if and only if $\be^{\rm T}\alpha=0$.
\end{lem}

\begin{proof}
    Let $\alpha=(a_1,\dots,a_n)^{\rm T}$, then we have
    \[\alpha^{\rm T}\alpha=\sum_{i=1}^na_i^2=(\sum_{i=1}^na_i)^2=(\be^{\rm T}\alpha)^2.\]
    So the Lemma follows immediately.
\end{proof}

\begin{lem}\label{lem:qii0}
    If the eigenvalue $\lambda_i\neq 0$, then we have $Q_{ii}=0$.
\end{lem}

\begin{proof}
By the construction of $P$ and $Q$, we see that the $(j,k)$ entry of $Q_{ii}$ is exactly $\alpha_{ij}^{\rm T}\alpha_{ik}$. Since $AP=PJ$, for $k=1,\dots,\ell_i-1$, we have $A\alpha_{i,k+1}=\lambda_i\alpha_{i,k+1}+\alpha_{ik}$, so for any $j,k=1,\dots,\ell_i-1$, we have
\[\alpha_{i,j+1}^{\rm T}\alpha_{i,k}=\alpha_{i,j+1}^{\rm T}(A-\lambda_i I)\alpha_{i,k+1}=\left((A-\lambda_i I)\alpha_{i,j+1}\right)^{\rm T}\alpha_{i,k+1}=\alpha_{i,j}^{\rm T}\alpha_{i,k+1}.\]
This implies that $Q_{ii}$ is of the form
\[\begin{pmatrix}
        a_1&a_2&\cdots&a_{\ell_i-1}&a_{\ell_i}\\
        a_2&a_3&\cdots&a_{\ell_i}&a'_{\ell_i-1}\\
        \vdots&\vdots&&\vdots&\vdots\\
        a_{\ell_i-1}&a_{\ell_i}&\cdots&a'_{3}&a'_{2}\\
        a_{\ell_i}&a'_{\ell_i-1}&\cdots&a'_{2}&a'_{1}
    \end{pmatrix},\]
To finish the proof of this lemma, it is enough to show that $\alpha_{ik}^{\rm T}\alpha_{ik}=0$ for any $k=1,\dots,\ell_i$ and $\alpha_{i,k+1}^{\rm T}\alpha_{ik}=0$ for any $k=1,\dots,\ell_i-1$.

For the first assertion, by Lemma \ref{lem:ealpha0}, it is enough to prove $\be^{\rm T}\alpha_{ik}=0$. we prove it by induction on $k$. When $k=1$, we have $A\alpha_{i1}=\lambda_i\alpha_{i1}$, so $\lambda_i \be^{\rm T}\alpha_{i1}=(A\be)^{\rm T}\alpha_{i1}=0$. This implies $\be^{\rm T}\alpha_{i1}=0$ since $\lambda_i\neq 0$. Now suppose the result is true for $k$, and we consider the case of $k+1$. Then by the induction hypothesis, we have
\[\lambda_i\be^{\rm T}\alpha_{i,k+1}=\be^{\rm T}A\alpha_{i,k+1}-\be^{\rm T}\alpha_{ik}=(A\be)^{\rm T}\alpha_{i,k+1}-\be^{\rm T}\alpha_{ik}=0.\]

For the second assertion, since $A\alpha_{i,k+1}=\lambda_i\alpha_{i,k+1}+\alpha_{ik}$, so together with the first assertion, we have
    \[\alpha_{i,k+1}^{\rm T}\alpha_{ik}=\alpha_{i,k+1}^{\rm T}A\alpha_{i,k+1}-\lambda_i\alpha_{i,k+1}^{\rm T}\alpha_{i,k+1}=0.\]
\end{proof}

\begin{lem}\label{lem:elli}
    If $\lambda_i=0$, then the entries of $Q_{ii}$ are $0$ except the position at $(\ell_i,\ell_i)$.
\end{lem}

\begin{proof}
With the same argument as the proof of Lemma \ref{lem:qii0}, we only need to prove that $\be^{\rm T}\alpha_{ik}=0$ and $\alpha_{i,k+1}^{\rm T}\alpha_{ik}=0$ for any $k=1,\dots,\ell_i-1$. For any $k<\ell_i$, we have $A\alpha_{i,k+1}=\alpha_{ik}$, so we get $\be^{\rm T}\alpha_{ik}=(A\be)^{\rm T}\alpha_{i,k+1}=0$. Moreover, we have $\alpha_{i,k+1}^{\rm T}\alpha_{ik}=\alpha_{i,k+1}^{\rm T}A\alpha_{i,k+1}=0$ since $A$ is symmetric and its diagonal entries are $0$.
\end{proof}

\begin{lem}\label{lem:q110}
    We have $Q_{11}=0$.
\end{lem}

\begin{proof}
    Note that $\alpha_{11}=\be$ by our assumption. By Lemma \ref{lem:elli}, the entry at the position $(1,\ell_1)$ is $0$ which is exactly $\alpha_{11}^{\rm T}\alpha_{1\ell_1}$. So we have $\be^{\rm T}\alpha_{1\ell_1}=\alpha_{11}^{\rm T}\alpha_{1\ell_1}=0$. Then by Lemma~\ref{lem:ealpha0}, we have $\alpha_{1\ell_1}^{\rm T}\alpha_{1,\ell_1}=0$. So the result follows from Lemma~\ref{lem:elli}.
\end{proof}

\begin{lem}\label{lem:qij0}
    If $\lambda_i\neq \lambda_j$, then $Q_{ij}=0$.
\end{lem}

\begin{proof}
    It is enough to prove that $\alpha_{ik}^{\rm T}\alpha_{jp}=0$ for all $k=1,\dots,\ell_i$ and $p=1,\dots,\ell_j$. We prove it by induction on $k+p$. When $k+p=2$, we have $k=p=1$, since $A\alpha_{i1}=\lambda_i\alpha_{i1}$ and $A\alpha_{j1}=\lambda_j\alpha_{j1}$, we have
    \[\lambda_i\alpha_{i1}^{\rm T}\alpha_{j1}=(A\alpha_{i1})^{\rm T}\alpha_{j1}=\alpha_{i1}^{\rm T}A\alpha_{j1}=\lambda_j\alpha_{i1}^{\rm T}\alpha_{j1}.\]
    This implies that $\alpha_{i1}^{\rm T}\alpha_{j1}=0$. Now we suppose that $k+p>2$. Without loss of generality, we may assume that $k>1$, then we have $A\alpha_{ik}=\lambda_i\alpha_{ik}+\alpha_{i,k-1}$. So by the induction hypothesis, we have
    \[\lambda_i\alpha_{ik}^{\rm T}\alpha_{jp}=\alpha_{ik}^TA\alpha_{jp}-\alpha_{i,k-1}^{\rm T}\alpha_{jp}=\alpha_{ik}^{\rm T}A\alpha_{jp}.\]
    If $p=1$, then $A\alpha_{j1}=\lambda_j\alpha_{j1}$, so we get $\alpha_{ik}^{\rm T}\alpha_{j1}=0$. If $p>1$, then by the induction hypothesis again, we get
    \[\lambda_i\alpha_{ik}^{\rm T}\alpha_{jp}=\alpha_{ik}^{\rm T}A\alpha_{jp}=\alpha_{ik}^{\rm T}(\lambda_j\alpha_{jp}+\alpha_{j,p-1})=\alpha_{ik}^{\rm T}\lambda_j\alpha_{jp}.\]
    This gives $\alpha_{ik}^{\rm T}\alpha_{jp}=0$.
\end{proof}

Now we are ready to give the proof of Theorem~\ref{thm:mainresult}.

\begin{proof}[Proof of Theorem \ref{thm:mainresult}]
    It suffices to prove that the order of each Jordan block $J_i$ does not exceed half of the algebraic multiplicity of the corresponding eigenvalue $\lambda_i$. By our construction, the order of $J_i$ is exactly the order of $Q_{ii}$. The main idea is the following obvious but key observation:

    \textbf{Observation:} For an invertible matrix of order $m$, the order of any zero submatrix does not exceed $m/2$.

    If $\lambda_i\neq 0$, without loss of generality, we may assume that the Jordan blocks $J_{s},\dots,J_{r}$ corresponding to the eigenvalue $\lambda_i$. Then the corresponding blocks of $Q$ is
    \[Q_i:=\begin{pmatrix}
        Q_{ss}&\cdots&Q_{rr}\\
        \vdots&&\vdots\\
        Q_{sr}&\cdots&Q_{rr}
    \end{pmatrix}.\]
    By Lemma \ref{lem:qij0}, we see that the matrix above is invertible. By Lemma \ref{lem:qii0}, we see that $Q_{ss},\dots,Q_{rr}$ are zero matrices. So by the claim above, we see that the order of $Q_{ii}$ does not exceed  $(\ell_s+\dots+\ell_r)/2$. This proves that the order of Jordan block $J_r,\dots,J_r$ does not exceed half of the algebraic multiplicity of the corresponding eigenvalue $\lambda_i$.

    Now we consider the eigenvalue $\lambda_i=0$. Without loss of generality, we may assume that the Jordan blocks $J_{1},\dots,J_{m}$ corresponding to the eigenvalue $0$ and the corresponding blocks of $Q$ is
    \[Q_1:=\begin{pmatrix}
        Q_{11}&\cdots&Q_{1m}\\
        \vdots&&\vdots\\
        Q_{m1}&\cdots&Q_{mm}
    \end{pmatrix}.\]
    In this case, we have no longer $Q_{ii}$ to be zero matrix except $Q_{11}=0$ by Lemma \ref{lem:q110}. So we need more works for $i\geq 2$. For any $i=2,\dots,m$, and $j=1,\dots,\ell_i-1$, by Lemma \ref{lem:elli}, the $(j,j)$ entry of $Q_{ii}$ is $0$ which is exactly $\alpha_{ij}^{\rm T}\alpha_{ij}$. So by Lemma \ref{lem:ealpha0}, we have $\be^{\rm T}\alpha_{ij}=0$. This implies that the position at $(1,\ell_1+\dots+\ell_{i-1}+j)$ is $0$, so together with Lemma \ref{lem:elli}, in the matrix $Q_1$, we find a zero submatrix of order $\ell_i$, i.e. the first row, the rows $\ell_1+\dots+\ell_{i-1}+1,\dots,\ell_1+\dots+\ell_{i-1}+\ell_{i}-1$ and the corresponding columns. So by the claim above, we see that $\ell_i\leq (\ell_1+\dots+\ell_m)/2$.

    Finally, let us consider the case that $n$ is odd. Let $\hat{A}:=\begin{pmatrix}
        A&0\\0&0
    \end{pmatrix}$. It is obvious to see that the characteristic polynomial of $\hat{A}$ is $x\phi(x)=\vp^2(x)$. Note that $\hat{A}$ has even order, we see that $\vp(\hat{A})=0$, this gives $\vp(A)=0$.

    This completes the proof.
\end{proof}

Actually, for odd $n$, we have a slightly stronger result in certain case.

\begin{thm}\label{odd}
    Let $A$ be an adjacency matrix of Eulerian graph of an odd order $n$. Let $\phi$ be the characteristic polynomial of $A$. Suppose that $\phi=x^k\phi_1^2$ over $\mbf_2$ with $k\geq 3$ and $2\nmid k$. Then $\varphi(x)/x=x^{\frac{k-1}2}\phi_1(x)$ is also the annihilating polynomial of $G$, i.e., $A^{\frac{k-1}{2}}\phi_1(A)=0$ over $\mbf_2$.

\end{thm}

\begin{proof}
The proof bears a resemblance to that of the even case, thus we will focus solely on elucidating the crucial distinctions between the two cases. In particular, our attention is directed towards the eigenvalue $\lambda_i=0$ since the proof is the same for nonzero eigenvalue. Assuming the same settings as outlined previously, a distinctive feature of the odd case is that we always encounter $\ell_1=1$ and $Q_{11}=1$. Otherwise, if $\ell_1>1$, then from $A\alpha_{12}=\alpha_{11}$ we have $\be^{\rm T}\be=\be^{\rm T}\alpha_{11}=(A\be)^{\rm T}\alpha_{12}=0$, but this is impossible since $\be^{\rm T}\be$ is odd. If $\ell_i> (k-1)/2$ for some $i$ where $k=\ell_1+\dots+\ell_m$ is the algebraic multiplicity of $0$, say $\ell_2$, then by Lemma \ref{lem:elli}, $Q_{22}$ has at most one non-zero entry at position $(\ell_2,\ell_2)$. Moreover, with the same argument as the proof of Theorem \ref{thm:mainresult}, the position at $(1,1+j),j=1,\dots,\ell_2-1$ of $Q_1$ are $0$. We see that $Q_1$ is of the form
    \[Q_1=\left(\begin{array}{c|ccccc|c}
        1&0&0&\cdots&0&\ast&\cdots\\
        \hline
        0&0&0&\cdots&0&0&\cdots\\
        0&0&0&\cdots&0&0&\cdots\\
        \vdots&\vdots&\vdots&&\vdots&\vdots&\\
        0&0&0&\cdots&0&0&\cdots\\
        \ast&0&0&\cdots&0&\ast&\cdots\\
        \hline
        \vdots&\vdots&\vdots&&\vdots&\vdots&\\
    \end{array}\right).\]
    So $Q_1$ has a $0$ submattrix of order $(\ell_2-1)\times (\ell_2+1)$, i.e. the rows $2,\dots,\ell_2$ and the columns $1,2,\dots,\ell_2+1$. However, such matrix $Q_1$ must be singular since $\ell_2>(k-1)/2$; a contradiction.
\end{proof}



\section*{Acknowledgments}
The research of the second author is supported by National Key Research and Development Program of China 2023YFA1010203 and National Natural Science Foundation of China (Grant No.\,12371357), and the third author is supported by Fundamental Research Funds for the Central Universities (Grant No.\,531118010622), National
Natural Science Foundation of China (Grant No.\,1240011979) and Hunan Provincial Natural Science Foundation of China (Grant No.\,2024JJ6120).

\end{document}